\documentclass{article}
\usepackage{times}
\input epsf

\newtheorem{Thm}{Theorem}

\newtheorem{Lem}{Lemma}

\newcommand{\vs}{\vspace{5 mm}}

\begin{document}

\title{Isonemal Prefabrics with Only Parallel Axes of Symmetry\thanks{Work on this material has been done at home and at Wolfson College, Oxford.
Will Gibson made it possible to draw the diagrams with surprising ease from exclusively keyboard input using Maple and Xfig. 
Richard Roth helped with the understanding of his papers and with this presentation.
An anonymous referee was helpful.
To them all I make grateful acknowledgement.}}
\author{R.~S.~D.~Thomas
}
\maketitle

\noindent St John's College and Department of Mathematics, University of Manitoba,

\noindent Winnipeg, Manitoba  R3T 2N2  Canada.
thomas@cc.umanitoba.ca

\noindent Other contact information: phone 204 488 1914, fax 204 474 7611.

\begin{abstract}
Isonemal weaving designs, introduced into mathematical literature by Gr\"un\-baum and Shephard, were classified into thirty-nine infinite sets and a small number of exceptions by Richard Roth. This paper refines Roth's taxonomy for the first ten of these families in order to solve three problems: which designs exist in various sizes, which prefabrics can be doubled and remain isonemal, and which can be halved and remain isonemal.
    \end{abstract}

\noindent AMS classification numbers: 52C20, 05B45.

\vs

\vs
\noindent  {\bf 1. Introduction}

\noindent Except for a short list of interesting exceptions, Richard Roth [1] has classified isonemal periodic prefabric designs into 39 infinite sets. 
The first ten have reflection or glide-reflection symmetries with parallel axes and no rotational symmetry. 
The remainder have reflection or glide-reflection symmetries with perpendicular axes, hence half-turns (22), or have quarter-turns (7). 
This paper is intended to reconsider this taxonomy for the first ten sets, refining it slightly to make it easier to use, and then to use it to examine three questions about prefabrics. 
This reconsideration presupposes the validity of Roth's taxonomy; the aim is to make it appear natural and to show it to be useful. 
The unanswered questions addressed are (for these designs) which prefabrics can be doubled and remain isonemal, which prefabrics can be halved and remain isonemal, and which orders there can be.
These questions are answered in the final three sections.
Answering depends on refining the taxonomy.

Since I am asking a reader to spend some time on a taxonomy, I ought to admit that no taxonomy that does not break its subject matter into individuals is not subject to further refinement.
Richard Roth used something slightly coarser than his taxonomy to study perfect and chromatic colouring of fabrics because he did not need its refinement for that purpose.
In this paper I am determining a finer taxonomy because it is needed in sections 7 and 8, although not in 6.
Section 7 actually requires a small further {\it ad hoc} refinement that I make nothing more of.

The argument of the paper is very simple.
It shows that knowing the symmetry groups of prefabric designs allows us to determine both facts about those designs (sections 6 and 7) and what designs there can be (section 8).
A lot of geometrical literature, including the paper [1] of Roth, stops short at crystallographic types of symmetry groups.
Useful as knowledge of types is---and it is crucial to the present work---it is the groups themselves that tell the story.
Roth's type determination is useful because it allows group determination with just a little additional effort.

As Roth observes beginning his subsequent paper [2] on perfect colourings, `[r]e\-cent mathematical work on the theory of woven fabrics' begins with [3], which remains the fundamental reference.
Roth's papers [1] and [2], however, contain the major advance from the fundamental work of Gr\"unbaum and Shephard.
In them he determines the various (layer, similar to crystallographic) symmetry-group types that periodic isonemal fabrics---actually prefabrics---can have and (in [2]) which of them can be perfectly coloured by striping warp and weft.
We are not concerned with striping warp and weft, but the other terms need to be defined.
A {\it prefabric,} as defined by Gr\"unbaum and Shephard [4], consists of two or more congruent layers (here only two) of parallel strands in the same plane $E$ together with a preferential ranking or ordering of the layers at every point of $E$ that does not lie on the boundary of a strand.
The points not on the boundary of any strand are naturally arranged into what are called here {\it cells}, in each of which one strand is uppermost.
The (parallel) strands of each layer are perpendicular to those of the other layer, making the cells square (Figure 1).
These square cells are taken here to be of unit area.
Mathematical literature on weaving has concerned exclusively {\it periodic} arrangements in the plane in the standard two-dimensional sense explained by Schattschneider in her exposition of plane symmetry groups [5].
There exists a non-unique finite region and two linearly independent translations such that the set of all images of the region, when acted upon by the group generated by these translations, reproduces the original configuration, which is assumed to be infinite in all directions for convenience.
Schattschneider gives the name {\it unit} to a smallest region of the plane having the property that the set of its images under this translation group covers the plane.
Such units are all of the same area, the {\it period}, but in general can be of a variety of shapes. 
Since our prefabric layers meet at right angles and the symmetry groups with which we shall be concerned here all have parallel axes of reflection or glide-reflection, the {\it lattice units,} that is, units whose vertices are all images of a single point under the action of the translation group, can be either rectangular or rhombic (Figure 3). 
Many of the rectangles have one set of parallel boundaries defined by the group but the perpendicular boundaries arbitrary in position, only the distance between them being dictated by the group.

The notion of symmetry group allows the definition of the term isonemal; a prefabric is said to be {\it isonemal} if its symmetry group is transitive on the strands, whose directions are conventionally chosen to be vertical, called {\it warps,} and horizontal, called {\it wefts.}
The distinction between prefabrics and fabrics can now be explained; a {\it fabric} is a prefabric that {\it hangs together}, that is, that does not {\it fall apart} in the sense that some warps or some wefts or some of each can be lifted off the remainder because they are not bound into a coherent network by the interleaving defined by the preferential ranking.

The standard way to represent the preferential ranking of the strands is to regard the plane $E$ as viewed from one side, from which viewpoint one or the other strand is visible in each cell.
By the {\it normal colouring} of warps dark and wefts pale, the visual appearance of the strands from a particular viewpoint becomes an easily understood code for which strand is uppermost.
An array of dark and pale congruent cells tessellating the plane is given a topological meaning, the {\it design} of the prefabric (Figure 1a).
\begin{figure}
\begin{tabbing}

\noindent
\epsffile{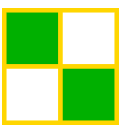}\hskip 10 pt\=
\epsffile{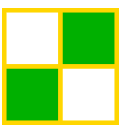}\hskip 10 pt\=
\epsffile{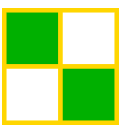}\hskip 10 pt\=
\epsffile{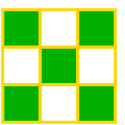}\hskip 10 pt\=
\epsffile{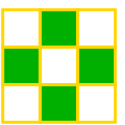}\\

\noindent (a)\>(b)\>(c)\>(d)\>(e)\\
\end{tabbing}

\noindent Figure 1. a. Fragment of plain weave. \hskip 10 pt b. Reverse of the fragment of Figure 1a as seen in a mirror. \hskip 10 pt c. Reverse of the fragment of Figure 1a as seen from behind. \hskip 10 pt d. Fragment of plain weave. \hskip 10 pt e. Reflection of the fragment of Figure 1d.
\end{figure}

If a finite region of the plane $E$ is viewed from behind by setting up a mirror beyond it and looking at $E$'s reflection, then what is seen is the strand perpendicular in direction and opposite or complementary in colour in corresponding positions as in Figure~1b compared with the front view of Figure~1a.
Figure 1b is in contrast to what one would see from the other side of $E$, which appears in Figure 1c, where the correspondence of cells of $E$ between Figure 1c and Figure 1a is obscured by the left-right reversal that causes mirror-image confusion in the real world.
For the sake of the correspondence, the reverse of a fabric will be represented as in Figure 1b; the mirror here being a simplifying device.
Another thing that the mirror does is to prevent symmetry axes from flipping between having positive slope and having negative slope.
As long as the strands are coloured normally, the reverse (so viewed) is just the colour-complement of the obverse and so is of no particular interest (`obverse' and `reverse' rather than `front' and `back' because of the arbitrariness of which is which, as of a coin).
Because one of the isometries that is used in weaving symmetries is reflection in the plane $E$, the reversal of which strand is uppermost at every non-boundary point, it is good to have an easy way to represent such reflection, denoted $\tau$.
As Figures 1a and 1b make clear, reversal of dark and pale represent the action $\tau$ adequately.
They also make clear why $\tau$ can be {\it in} a symmetry but cannot itself {\it be} a symmetry.

If Figure 1d is reflected in either diagonal, the diagram remains Figure 1d, but, because warps become wefts under such a reflection and vice versa, Figure 1d is not the design of the reflected fabric.
The representation of the reflected fabric is Figure 1e, correctly indicating that the reflection is not a symmetry of the fabric.
If, however, $\tau$ is combined with the operation, then the colours are reversed again, Figure 1d is restored, and the combination of (two-dimensional) diagonal reflection ({\it within} the plane) and $\tau$ (three-dimensional reflection {\it in} the plane) is seen to be a symmetry of the fabric because it preserves its design.

The reference numbers of catalogued fabrics and prefabrics require explanation before use. 
Because the prefabrics of interest are periodic in a two-dimensional sense, they are also periodic along each strand. 
Because they are isonemal, they have the same period along every strand, which is called here {\it order} [3] to prevent confusion with two-dimensional period, which often differs numerically.%
\footnote{An order length of one strand is a unit but not a lattice unit for prefabrics of genus I, but that length of a pair of adjacent strands is needed to compose a unit for pure genus II, 
making two-dimensional period either one or two 
times order respectively.
The categorization of these designs into genera will be explained shortly.}
Fabrics were catalogued by Gr\"unbaum and Shephard, sorted first under order, and then by binary index, which is the sequence of pale and dark cells of order length represented by 0s and 1s and chosen to be a minimum, then given an arbitrary sequence number.
The first catalogue [6] listed fabrics other than twills up to order 8 (plus 13) and an extension [7] likewise up to order 12 (plus 15 and 17).

Prefabrics that fall apart%
\footnote{The subject of several investigations: Clapham [8] and Ens [9] for prefabrics in general, W.D.~Hoskins and Thomas [10] and Gr\"unbaum and Shephard [6] for isonemal prefabrics, and Delaney [11] for non-isonemal prefabrics.}
were catalogued by J.A.~Hoskins and Thomas [12] up to order 16---numerically the same way but with an asterisk to indicate that they are not fabrics.
The omission of {\it twills} from catalogues, despite their all being isonemal, stems from their extreme simplicity; the over-and-under sequence of each strand is that of its neighbour, say below, shifted by one place always in the same direction, e.g., Figure 5b.
This row-to-row shift is called the {\it offset.}
Because twills do not appear in the catalogues, a twill is referred to by its over-and-under sequence along a strand as in the caption to Figure 5b.

The designs to be discussed here include twills and the simplest extension of the twills, called twillins [3].
{\it Twillins} have each row shifted more than one place always in the same direction up the design ({\it offset} greater than 1).
If from row to next row, in addition to the shift, the colours are complemented, the design was at first called a `colour-alternate twillin' but now is referred to as being of {\it genus II} [6].
{\it Genus I} contains all twills and twillins without colour reversal.
These two genera (and the other three not discussed here) can overlap: if an isonemal design has pale and dark cells balanced in number along a strand (and therefore along all strands), colour reversal need not change the sequence (e.g., plain weave in Figure 1).
A design of genus II that is not also of genus I is called, as in footnote 1, of {\it pure} genus II.

Design is a mathematical model of woven fabrics. 
It may be helpful to consider the comparison between a fabric and the models briefly.
Consider a fabric fragment with design illustrated in Figure 2a (12-35-1 with some mirrors and boundaries of a lattice unit indicated) thought of for a moment as fixed in a horizontal plane and viewed from above and from the south.
\begin{figure}
\begin{tabbing}
\noindent\epsffile{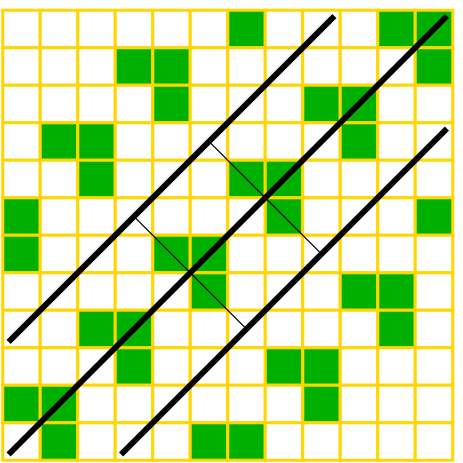}\hskip 10.5 pt\=\epsffile{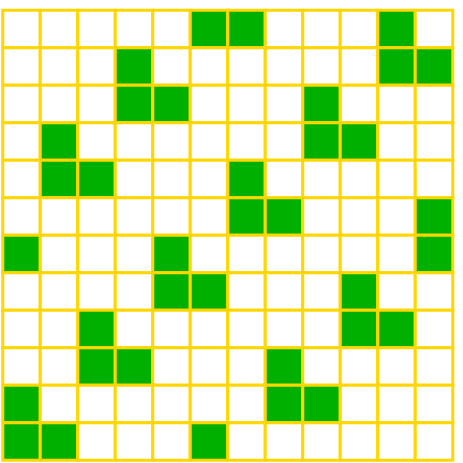}\\

\noindent (a)\>(b)\\

\noindent\epsffile{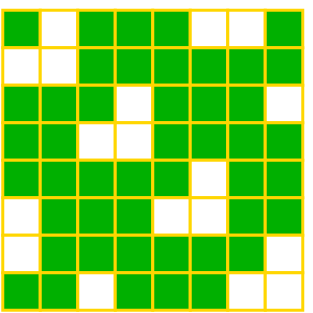}\hskip 5 pt\=\epsffile{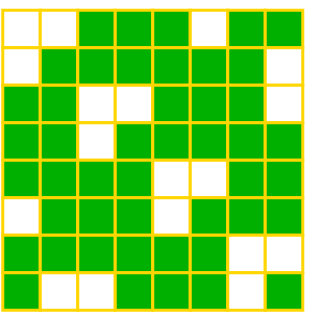}\hskip 5 pt\=\epsffile{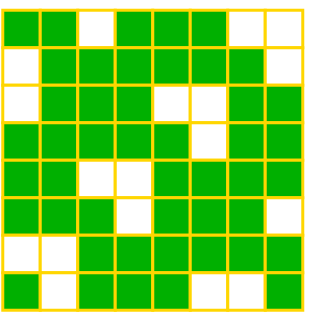}\\  

\noindent (c)\>(d)\>(e)\\

\noindent\epsffile{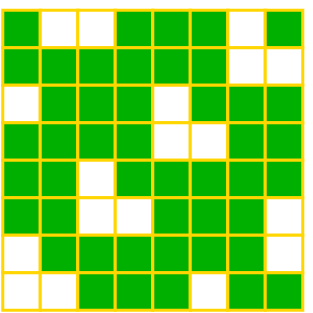}\>\epsffile{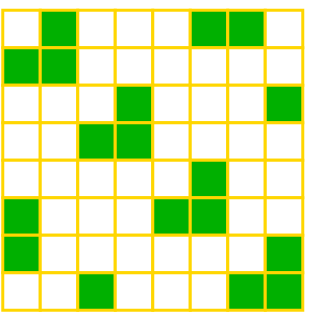}\>\epsffile{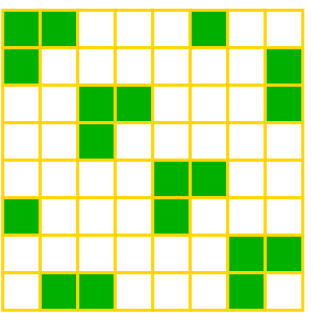}\\

\noindent (f)\>(g)\>(h)

\end{tabbing}
\noindent Figure 2. Fragment of fabric 12-35-1, stationary as viewer takes up viewpoints. a. View from above and the south. \hskip 10 pt b. View from above and the north. \hskip 10 pt c. View from above and the east. \hskip 10 pt d. View from above and the west. \hskip 16 pt e. View from below and south. \hskip 16 pt f. View from below and north. \hskip 10 pt g. View from below and east. \hskip 10 pt h. View from below and west.
\end{figure}
In the terms of the previous paragraph, it is a twillin with offset 5, order 12, and genus I.
Because the fabric does not have half-turn symmetry (which would make its genus III, IV, or V), viewed from the north it looks quite different (Figure 2b). 
Looked at from the east and from the west (Figures 2c and 2d respectively), it looks different in a different way because of the colouring convention (which is relative to the view).
The four views from below (and south, north, east, and west) are different again, complementary to the first four (Figures 2e--h, respectively).
The four views from below would of course also be seen if the fabric fragment were turned over and still viewed from above, but which diagram is which would depend on how the fragment was turned over.
There are eight different designs illustrated here all representing the same piece of physical fabric, a situation somewhat opposed to the usual simplification in mathematical models.
For no good reason, illustrations here will have axes of positive slope. 
There is a handedness to all of these designs that the physical fabrics modelled lack. 
Whether the designs illustrated here are right-handed or left-handed is rightly left to the reader.

In the trivial prefabric, all the wefts pass over all the warps or vice versa.
The lattice unit is a single cell, conventionally pale or dark respectively.
As Roth points out `for completeness', its symmetry group is a layer group that is a subgroup of the direct product of a two-dimensional group of type $p4m$ (having several mirror symmetries in addition to rotational symmetry) with this lattice unit and the group in the third dimension $\{e,\tau\}$, where $e$ is the identity.
For this simplest prefabric, all interchanges of warp and weft require $\tau$ to make them symmetries.
So the symmetry group consists of operations $(g,e)$, where $g$ does not interchange warp and weft, and $(g,\tau)$, where $g$ involves a quarter-turn or reflection or glide-reflection with axis at $45^\circ$ to the strand directions.
All other prefabrics have lattice units that are not single cells.
The full symmetry group is a three-dimensional layer group $G$.
Associated with $G$ is a planar group algebraically isomorphic 
to $G$ consisting of all $G$'s elements $g$ whether paired with 
$e$ or $\tau$, the planar projection of $G$, denoted $G_1$ by Roth. 
$G$'s {\it side-preserving subgroup} $H$, consisting of the 
elements $(g, e)$ and omitting $(g, \tau)$ elements is also 
isomorphic%
\footnote{This stronger geometrical sense of isomorphism is explained by 
Gr\"unbaum and Shephard [13], p.~38, and can be called {\it geometrical isomorphism}.}
 to its planar projection 
$H_1$ consisting of those elements $g$ paired with $e$ in $H$.
$H_1$ will be referred to as the side-preserving subgroup of $G_1$.
Figure 2a illustrates the common lattice unit for $G_1$ and $H_1$ for the fabric 12-35-1 and axes of reflective symmetry that are in $G_1$ but not in $H_1$.
Note that, while the lattice unit for this $H_1$ can be chosen to be the same as that for $G_1$, it is often larger, corresponding as it does to a subgroup of $G_1$.

For the prefabrics with only parallel symmetry axes, the planar projection $G_1$, the group comprising the set of $g$s, is of crystallographic type $pg$, $pm$, or $cm$, that is, it has parallel glide-reflection axes, or axes of reflection, or both (alternating) respectively.
In order for a glide-reflection to be a symmetry of the prefabric, it may or may not have to be combined with reversal of the sides of the prefabric, $\tau$. 
Mirror symmetry must always be combined with $\tau$ because any cell through which the mirror passes must be its own image in the symmetry but the reflection alone, reversing warp and weft, would reverse its conventional colour. So $\tau$ is needed to restore it.
This means that there is never mirror symmetry in $H_1$. 
It is $H_1$ that determines the two-dimensional period under translation alone.
For these prefabrics, $H_1$ is of type either $p1$, which is a group of translations only, or $pg$.
\begin{figure}
\noindent
\epsffile{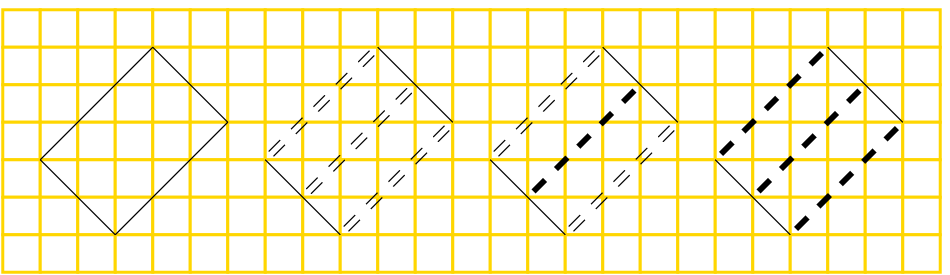}

\noindent \hskip 0.16 in (a)\hskip 0.73 in (b)\hskip 0.73 in (c)\hskip 0.73 in (d)

\noindent
\epsffile{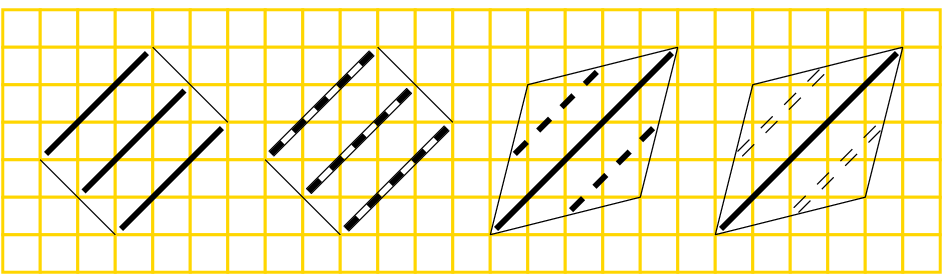}

\noindent \hskip 0.16 in (e)\hskip 0.73 in (f)\hskip 0.73 in (g)\hskip 0.73 in (h)

\noindent Figure 3. Representations of side-preserving-subgroup and symmetry-group types. a. Type $p1$ (side-preserving subgroup only). \hskip 10 pt b, c, d. Type $pg$ (three versions: $pg/-$, $pg/pg$, and $pg/p1$). \hskip 10 pt e, f. Type $pm$ (two versions: $pm/p1$ and $pm/pg$). \hskip 10 pt g, h. Type $cm$ (two versions: $cm/p1$ and $cm/pg$).
\end{figure}
The possible side-preserving subgroups just mentioned and their lattice units, dimensions aside, are illustrated in Figures 3a and b, where glide-reflection axes without $\tau$ are hollow dashed black lines.
Later in Figure 3, dashed axes of glide-reflections with $\tau$ are filled not hollow (Figures 3c, d, and g) like the undashed mirrors (Figures 3e, g, and h). 
The dashed filling of hollow lines in Figure 3f indicates a mirror (with $\tau$) that is also an axis of glide-reflection without $\tau$.
Thin black lines are just boundaries of lattice units, outlining or completing the outline of rectangles or rhombs.
Symbols in the caption of Figure 3 will be explained in the next section.
The rectangular and rhombic shapes of the lattice units have a maximum dimension in the direction of the axes.
It will be referred to as the unit's {\it length}, and the maximum dimension perpendicular to the axes will be referred to as its {\it width} even when width is greater than length and when, in the rhombs, width---like the length---is a diagonal.
\vs

\noindent {\bf 2. Symmetry operations}

\noindent Roth has determined the types of the layer groups that are the three-dimensional symmetry groups of prefabrics in terms of the standard crystallographic types of the layer groups' planar projections.
Corresponding to these types of layer groups and their corresponding two-dimensional groups are what I am calling species of prefabric.
Something to note about the groups is that they are unavoidably subgroups of other groups, but for classification purposes it is helpful to have the species not overlap.
A {\it species} of isonemal prefabrics is the set of prefabrics with symmetry groups of a single Roth type or subtype.
The Roth symmetry-group types and therefore the species of prefabrics here are numbered 1 to 10; some will be subdivided, as 7 into $7_o$ and $7_e$, where the subdivisions will be indicated by a subscript.
A {\it family} of isonemal prefabrics is the set of prefabrics with the same symmetry group, not group type.
Families are obviously finite; all the species except the few interesting exceptions to which I alluded in the first sentence contain infinitely many families.
The exceptions all have order 4 or less.
As such prefabrics are exceptions to many otherwise general statements, attention here is limited to orders greater than 4.

The way I have chosen to develop Roth's categorizing of prefabrics with only parallel axes of symmetry is to consider $H_1$ first. 
There are only two $H_1$ diagrams except for dimensions (Figures 3a and 3b).
The simplest symmetry-group division divides the prefabrics into species 2, 5, and 8 having $H_1$ of type $p1$ and the others having $H_1$ of type $pg$.
What of $G_1$?
It is possible that $G_1=H_1$, but only if $H_1$ is of type $pg$.
This gives the $G_1/H_1$ type combination%
\footnote{The notation is a Coxeter notation specifying the layer group $G$, not an indication of a factor group.}
that Roth denotes $pg/-$, the dash meaning that $H_1$ is not just of the same type as $G_1$ but is the same group.
The combination characterizes Roth type 1.

To either $p1$ or $pg$ can be added glide-reflections with $\tau$ giving groups of type $pg$ of either the second sort illustrated in Figure 3d ($pg/p1$, Roth type 2) or of the third (mixed) sort illustrated in Figure 3c ($pg/pg$, Roth type 3 or 4 depending on further considerations).
The repetition of $pg$ indicates that, while the planar projections $G_1$ and $H_1$ are of type $pg$, they are not the same group.
The lattice unit of groups of the $p1$ and $pg$ types can be used for groups of type $pm$ depending on whether, in the side-preserving subgroup, the mirrors simply disappear ($pm/p1$, Roth type 5, Figure 3e) or become axes of side-preserving  (no $\tau$) glide-reflection ($pm/pg$, Roth type 6 or 7 depending on further considerations, Figure 3f).
One produces the first sort of $cm$ in Figure 3g, by adding to a rhombic lattice unit for $p1$ (such lattice units need only be parallelograms, but here they are rectangular except in this case) both mirrors and glide-reflection axes with $\tau$ ($cm/p1$, Roth type 8).
One produces the second sort of $cm$ in Figure 3h by adding to type-$pg$ $H_1$ of Figure 3b mirrors between the axes of side-preserving glide-reflection.
This process gives $cm/pg$ with $G_1$ lattice unit of a quite different shape but the same length, Roth type 9 or 10 depending on further considerations. (See Figure 7b for both lattice units in one diagram.)
The further considerations now mentioned three times have to do with the placement of the axes with respect to the grid of cells, the subject to which we must now turn.

A mirror can intersect the boundary of a cell only at corners, either singly or in diagonally opposite pairs, transforming the cell into either a diagonally adjacent cell or itself. 
A symmetry group containing both reflection in an axis and a
translation along that axis contains, as their composition in
either order, a glide-reflection with the same axis and glide
the length of the translation.
This possibility illustrates that glide-reflection axes may have what I call {\it mirror position}.
But as well, such axes can be displaced to intersect cell boundaries at the mid-points of adjacent sides and pass through the right part of a cell, then the left part of a cell, and so on alternately.
This extra freedom in the placement of glide-reflection axes leads to Roth's types 4 and 10 rather than 3 and 9, and it is the source of the refinement of his taxonomy proposed here.
Roth's example of a species-1 fabric illustrates the mirror-position possibility (Figure 4a), but Figure 4b illustrates the other position.
\begin{figure}
\noindent
\epsffile{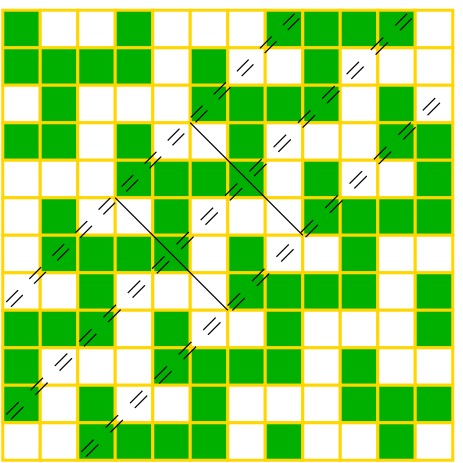}\hskip 10 pt\epsffile{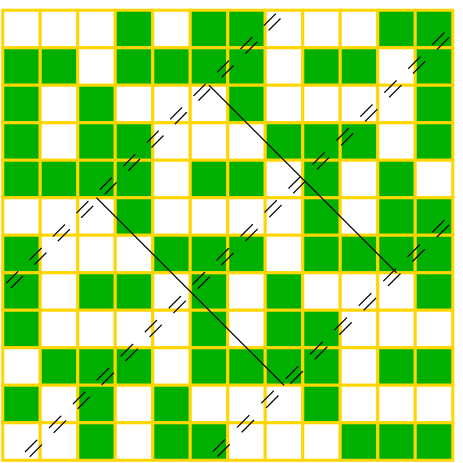}

\noindent (a)\hskip 1.8 in (b)

\noindent Figure 4. $G_1$ (same as $H_1$) lattice units of
fabrics of species 1.\hskip 10 pt a. 12-183-1 with glide-reflection axes in mirror position. \hskip 10 pt b. Order-30 fabric with glide-reflection axes not in mirror position.
\end{figure}
Figure 4b illustrates that a glide-reflection with axis not in mirror position reflects axes of other glide-reflections not in mirror position among themselves, as do Figures 8b, 10, 12b, and 13.

\begin{Lem}{Glide-reflection axes of a prefabric with symmetry group of crystallogra\-phic type $pg$ are all in mirror position or all not in mirror position.}
\end{Lem}

\noindent {\it Proof.} To see that a glide-reflection with axis in mirror position cannot reflect axes not in mirror position among themselves, it suffices to consider exactly what it would have to do.
For definiteness, consider two glide-reflection axes not in mirror position bounding a lattice unit of a symmetry group of type $pg$ as in Figure 4b.
Because each axis is translated perpendicularly to the other axis (as well as in many oblique directions), they both have corresponding parallel paths through the right part of a cell, then the left part of a cell, then the right part of a cell, and so on, in either direction.
The glide-reflection with axis central to the lattice unit must transform each axis to the other one.
It must therefore be half-way between them, and this location is not a mirror position but is a similar or opposite (depending on the distance between the other axes) path through the right part of a cell, then the left part of a cell, and so on.
The lemma is proved.

Glide-reflection axes can be displaced from mirror position when they are combined with mirrors because there is not, in such groups, a translation from one axis to the next that is perpendicular to the axes.
Such perpendicular translations as there are relate {\it alternate} glide-reflection axes as in the fabrics in Figures 5, 7, 9, 16b and 17; the translations between adjacent glide-reflection axes are exclusively oblique (the sides of the rhombs indicating the least oblique).
\begin{figure}
\noindent
\epsffile{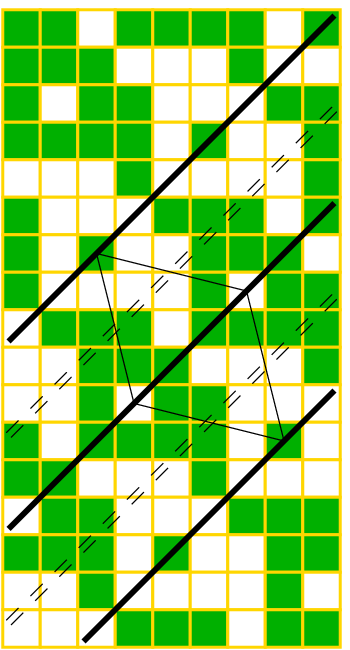}\hskip 10 pt\epsffile{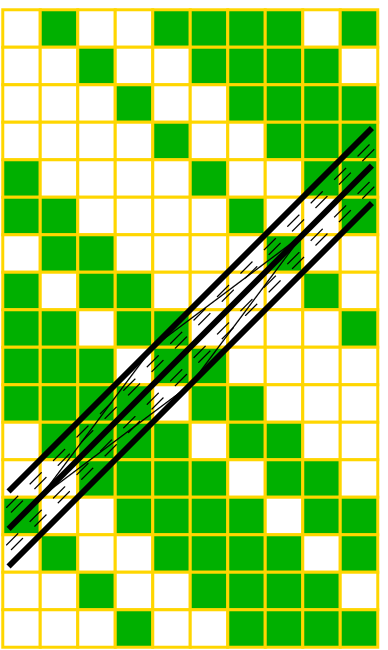}

\noindent (a)\hskip 1.35 in (b)

\noindent Figure 5. Species-10 fabrics with $G_1$ lattice units.\hskip 10 pt a. Roth's order-30 example of a species-10 fabric. Figure 5 of [1]. \hskip 10 pt b. The 1/2/4/1/2/4 twill.
\end{figure}

\vs

\noindent {\bf 3. Lattice units} 

\noindent Since lattice units have already been marked on figures, it would be as well to stress that the location of unit boundaries other than axes has no geometrical significance except for how far apart they are.
When one boundary is chosen, another is its closest parallel image under the relevant group, $G_1$ or $H_1$.

With the sole exception of species 5 (of type $pm/p1$), which has only mirror symmetry (including $\tau$) besides translations in $G_1$, and so only translations in $H_1$, the lengths of lattice units are determined by the glides of the glide-reflections present in $G_1$ of all prefabrics considered here.
When a glide-reflection axis is in mirror position, the glide must be an integer multiple of the cell diagonal $\delta$ to transform cells to cells (e.g., along the axis).
The minimum translation along the axis that is a symmetry, with or without $\tau$, is accordingly twice the glide and so an even multiple of $\delta$.
The length of a lattice unit for such glide-reflection, whether rhombic or rectangular, is accordingly an even multiple of $\delta$.
On the other hand, when a glide-reflection axis is not in mirror position, the relevant increment is $\beta \equiv \delta/2$.
Along the axis (and hence everywhere), cells can be transferred to cells if the glide is an integer multiple of $\beta$.
But to be a glide-{\it reflection,} the multiple must be odd.
The length of a lattice unit in such a group must therefore be twice an odd multiple of $\beta$, which is an odd multiple of $\delta$.
We have proved the following lemma.

\begin{Lem} According to whether the axis of a glide-reflection of an isonemal prefabric design is in mirror position or not, the length of the $H_1$ lattice unit is an even or odd multiple of the cell diagonal.
\end{Lem}

We can see now that the lengths of glides determine the lengths of lattice units for a group including a glide-reflection and for its side-preserving subgroup, which includes the square of any glide-reflection in the group.
For prefabrics of species 6 and 7, whose mirrors become axes of side-preserving glide-reflection in their side-preserving subgroups, it should be noted {\it how} the glide determines the length of the $G_1$ lattice unit.
Because the axis of the glide-reflection in $H_1$, which advances cells of the diagram by the glide and allows their colours to be complemented by the colouring convention (no $\tau$), is coincident with the axis of the diagram's mirror symmetry (not in $H_1$), the glide-reflection effects what looks like a translation of the diagram.
Such a translation with colour complementation is not a symmetry of the diagram because of the colour complementation; the corresponding symmetry is what Roth calls a $\tau${\it -translation.}
\begin{figure}
\noindent
\epsffile{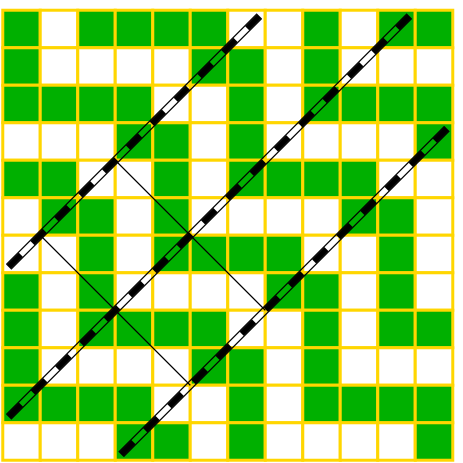}\hskip 10 pt\epsffile{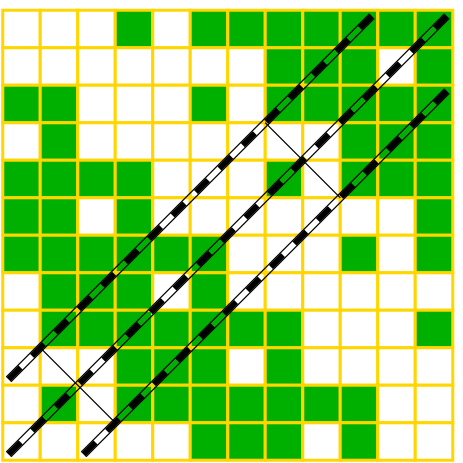}

\noindent (a)\hskip 1.8 in (b)

\noindent Figure 6. Roth's species-6 examples with $H_1$ lattice units. The $G_1$ lattice units have half the length.\hskip 10 pt a. 8-11-1; $\ell =1$, $w=4$. \hskip 10 pt b. 12-23-1; $\ell =3$, $w=2$.
\end{figure}
It is illustrated by Roth's species-6 examples in Figure 6, where one half of the $H_1$ lattice unit $\tau$-translates along the axes to the other half.
{\it Cf.} Figures 15 and 16a.

Something related occurs in species 9 (Figure 7) and 10 (Figure 5), where again there are side-preserving glide-reflec\-tions in $H_1$ and mirrors in $G_1$, but they are not coincident.
\begin{figure}
\noindent
\epsffile{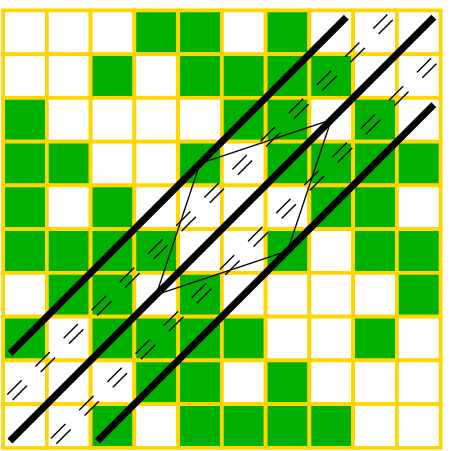}\hskip 10 pt\epsffile{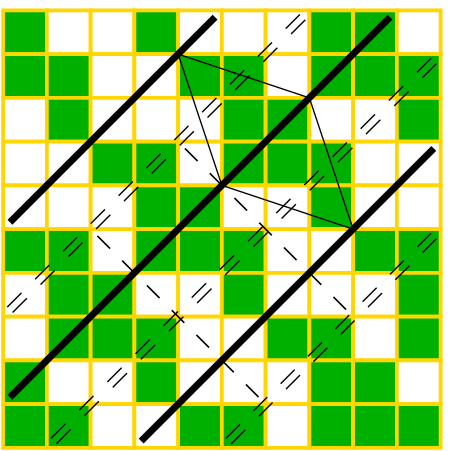}

\noindent (a)\hskip 1.74 in (b)

\noindent Figure 7. Roth's species-9 examples with $G_1$ lattice units.\hskip 10 pt a. 8-11-2. \hskip 10 pt b. 8-19-2 ($H_1$ lattice unit dashed).
\end{figure}
Mirrors are interchanged by the glide-reflection half-way between them.
Because of the mirror symmetry, the colour-reversing glide-reflection effects what looks like a translation but oblique to the direction of the axes because of the non-coincidence of the mirror and the glide-reflection's axis. 
This effect explains the $\tau$-translations along the boundaries of the rhombic lattice unit of these types and so how the rectangular lattice unit of $H_1$ has twice the area although the same length as the rhombic lattice unit of $G_1$.

We have seen two ways in which the lattice unit of $H_1$ can be larger than that of $G_1$. To see a third, we must change our attention to the width of the lattice unit of $G_1$ of types 3 and 4, formed by adding, between each adjacent pair of axes of side-preserving glide-reflection as shown in Figure 3b, an axis of side-reversing glide-reflection as shown in Figure 3c.
The $H_1$ lattice unit, as illustrated in Figure 8, can be taken to comprise two $G_1$ lattice units side by side.
\begin{figure}
\noindent
\epsffile{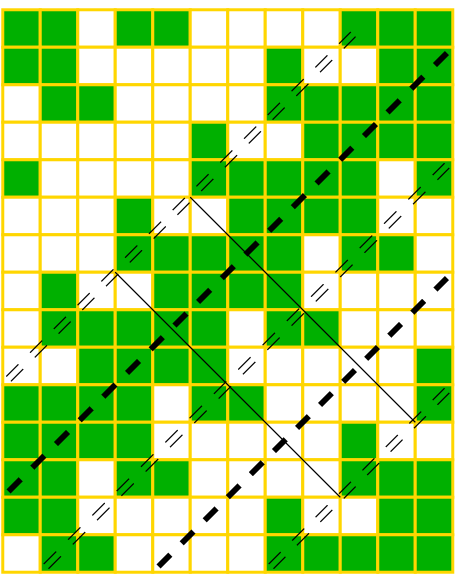}\hskip 10 pt\epsffile{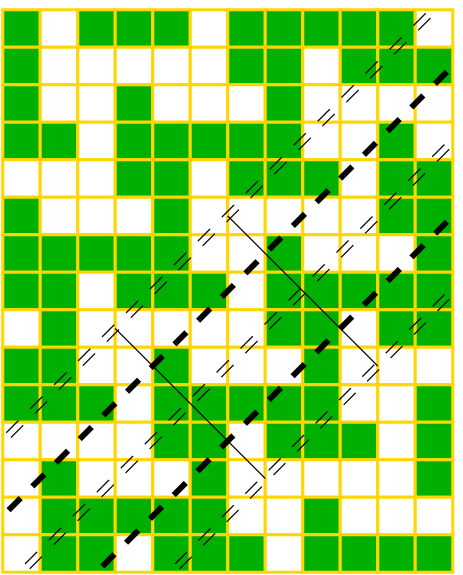}

\noindent (a)\hskip 1.79 in (b)

\noindent Figure 8. Roth's $pg/pg$ examples with side-by-side $G_1$ lattice units composing $H_1$ lattice units. 
\hskip 10 pt a. Species 3: 12-79-1. \hskip 10 pt b. Species 4: 24-282555. Figure 4 of [1].
\end{figure}
This convention is adopted in order to have side-preserving glide-reflection axes as the sides of both $G_1$ and $H_1$ lattice units as they are illustrated in Figures 3b and c.
These side-by-side $G_1$ lattice units are not related by translation as they would be if they were $H_1$ lattice units but rather by $\tau$-translations.
Figure 8 shows designs of species 3 (axes in mirror position) and 4 (axes not in mirror position).

While the length $\ell\delta$ of a rhombic lattice unit is twice the length of the glides of the symmetry group and the width $w\delta$ is twice the distance between adjacent mirrors, only certain combinations of $\ell$ and $w$ are feasible.
The constraints on mirror position and glide lengths require only that $\ell$ and $w$ be integers.
Without loss of generality, let the centre of a lattice unit, which lies on a mirror, lie at a cell corner as illustrated in Figure 7b.
(Any point on the mirror is a possible choice, e.g., a cell centre as in Figure 7a.)
The distance $w\beta$ from the unit centre to the adjacent mirror can be either even as in Figure 9a or odd as in Figure 9b.
If $w$ is even, then the corners of the lattice unit on the adjacent mirrors fall at cell corners; if $w$ is odd, then the corresponding corners fall at cell centres.
The other corners of the lattice unit must also fall at cell corners or cell centres respectively.
This means that $\ell\beta$ along the central mirror must have $\ell$ even when $w$ is even and $\ell$ odd when $w$ is odd.
The feasibility constraint on $\ell$ and $w$ is that their parities be the same.
\begin{figure}
\noindent
\epsffile{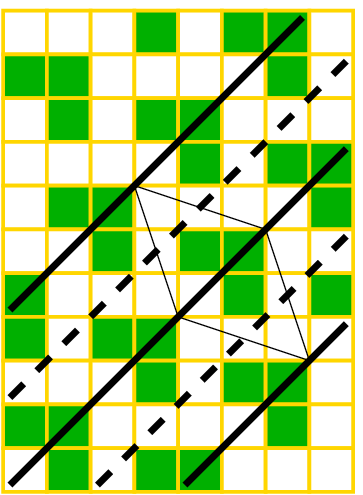}\hskip 10 pt\epsffile{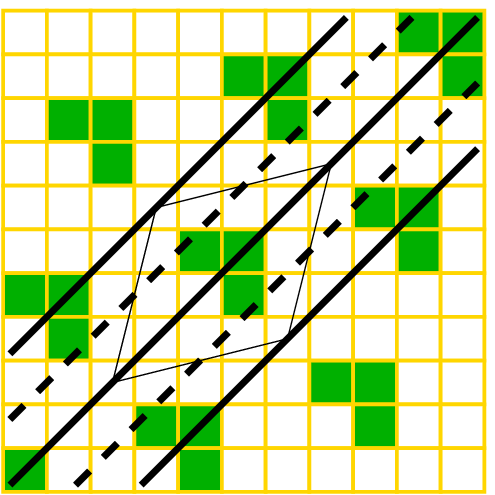}

\noindent (a)\hskip 1.4 in (b)

\noindent Figure 9. Roth's examples of species 8 with $G_1$ and $H_1$ lattice units. a. 8-11-6 with $(\ell , w)=2$. \hskip 10 pt b. 15-35-1 with $(\ell , w)=1$.
\end{figure}
\vs

\noindent {\bf 4. Isonemality} 

\noindent For a symmetry group to make a prefabric isonemal, a corner of its lattice unit, if chosen to be in a cell, must have translates in every other strand. 
The side-preserving subgroup $H_1$ need not and does not always do this.

Let rectangular $G_1$ lattice units be measured in the diagonal $\delta$ of a cell.
A unit $\ell\delta$ long by $w\delta$ wide will have a corner image in every strand when the greatest common factor of $\ell$ and $w$ is 1, equivalently $\ell$ and $w$ are relatively prime or there exist integers $c$ and $d$ such that $c\ell + dw =1$.
The last condition is a kind of recipe for how many iterated translations $|c|$ and $|d|$ in perpendicular directions are required to move a corner to an adjacent parallel strand. 
A dimension can sometimes be $\delta$ as in Figure 5b. 
In contrast, the similar parameters for $H_1$ can have greatest common divisor 2 and often do (Figures 6, 7, 8a, and 13).
This can be summed up in a lemma.
\smallskip

\begin{Lem} A necessary and sufficient condition that a prefabric of species $1$ to $7$ with a rectangular $G_1$ lattice unit $\ell\delta$ by $w\delta$ be isonemal is that $(\ell , w)=1$.
\end{Lem}

If rhombic $G_1$ lattice units that tessellate the plane without colour reversal ($cm/p1$, Roth type 8) or with colour reversal ($cm/pg$, Roth types 9 and 10) are measured similarly along their diagonals, then the necessary condition for isonemality may be $(\ell , w)$ is either 1 or 2.
Roth gives examples of both for species 8, illustrated in Figure 9.
We also have a lemma for these other species.

\begin{Lem} A necessary and sufficient condition that a prefabric of species $8$ to $10$ with a rhombic $G_1$ lattice unit of length $\ell\delta$ and width $w\delta$ be isonemal is either that $\ell$ and $w$ be odd and $(\ell , w)=1$ or that $\ell$ and $w$ be even, $(\ell /2, w/2)=1$, and $\ell /2$ and $w/2$ differ in parity.
\end{Lem}

\noindent{\it Proof.} 
If $\ell\delta$ and $w\delta$ are the lengths of rhombic diagonals, then the diagonals can be thought of as vectors $(\ell , \ell)$ and $(w,-w)$.
The sides of the rhombs can be thought of as vectors $\textstyle{\left(\frac{\ell +w}{2},\frac{\ell -w}{2}\right)}$ and  $\textstyle{\left(\frac{\ell -w}{2},\frac{\ell +w}{2}\right)}$.
The feasibility constraint of section 3 is that $\ell$ and $w$ must have the same parity.

For isonemality it is necessary and (with the reflective symmetry being supposed) sufficient that horizontally and vertically a rhomb corner in one strand be translatable to an adjacent strand as for Lemma 3.
Here this condition is that integers $e$ and $f$ exist such that $\textstyle{e\frac{\ell +w}{2}+f\frac{\ell -w}{2}=1}$.
For $\ell $ and $w$ both odd, $\textstyle{e\frac{\ell +w}{2}+f\frac{\ell -w}{2}=1}$ is just a rearrangement of  $\textstyle{\frac{e+f}{2}\ell+\frac{e-f}{2}w=1}$.
So for $\ell $ and $w$ odd, isonemality is equivalent to the relative primality of $\ell $ and $w$.
In Figure 9b, $\ell=5, w=3$, and so $e=1, f=-3$.
When $\ell $ and $w$ are even, $c\ell  + dw =1$ is impossible, and so a different condition is needed.
Since another rearrangement of $\textstyle{e\frac{\ell +w}{2}+f\frac{\ell -w}{2}=1}$ is 
$\textstyle{(e+f)\frac{\ell}{2}+(e-f)\frac{w}{2}=1}$, divisions possible because $\ell$ and $w$ are both even, the condition required includes the relative primality of $\ell /2$ and $w/2$.
In Figure 9a, $\ell=2, w = 4$, and so $e=1, f=2$.
The specific form of the coefficients $e\pm f$, themselves unavoidably of the same parity whether $e, f$, were odd-odd, odd-even, or even-even, forces $\ell /2$ and $w/2$ to {\it differ in parity as well as being relatively prime,} since if $\ell /2$, $w/2$, as well as $e+f$, $e-f$, were of the same parity, the displayed equation could not be satisfied.
(Since $e\pm f$ must both be odd, $e$ and $f$ will also differ in parity.)

\vs

\noindent {\bf 5. Symmetry groups} 

\noindent The preceding observations make it possible
to determine each of the infinite families of symmetry groups $G$ that isonemal prefabrics with no rotational symmetry have in terms of the length and width of their lattice units.
We shall work through Figure 3 from 3b to 3h.

\noindent {\it 1.}\hskip 8 pt When $G_1=H_1$ is of type $pg$, the axes of side-preserving glide-reflection can be in mirror position or not.
If in mirror position, then the lattice unit is an even multiple $\ell\delta$ long.
By the standard isonemality constraint, the width must be $w\delta$ with $(\ell,w)=1$ and so $w$ odd.
As Roth remarks [7, p.~318], $w>1$ or a twill with mirror symmetry would result, contradicting the type of the symmetry group.
This is species $1_m$ (m for mirror position).
Roth's species-1 examples have symmetry of this subtype; a twillin with offset 5 is illustrated in Figure 4a. 
Period and order are $2\ell w$, always divisible by 4.
Genus is I.

If the glide-reflection axis is displaced, then the length of the lattice unit must be an odd multiple $\ell\delta$.
For isonemality, the width must be $w\delta$ with $(\ell,w)=1$, allowing $w$ to be even or odd.
There are accordingly two further species $1_e$ and $1_o$, respectively, with parameters as follows: For $1_e$, parameters are as for $1_m$ except for the omission of $\ell=3$, $w=2$, which forces the lattice unit to have more symmetry, e.g., 12-17-1 with group of Roth type 16.
This is analogous to barring $w=1$ above. 
Period and order are $2\ell w$, always divisible by 4.
Genus is I.

For $1_o$, Figure 4b showed the example $\ell =3$, $w=5$. 
Figure 10 shows $\ell =3$, $w=4$ and $\ell =5$, $w=3$.
Parameters $\ell =1$, $w=2$ and $\ell =3$, $w=2$ are too small.
Period and order are $2\ell w$, always even.
Genus is I.
\begin{figure}
\noindent
\epsffile{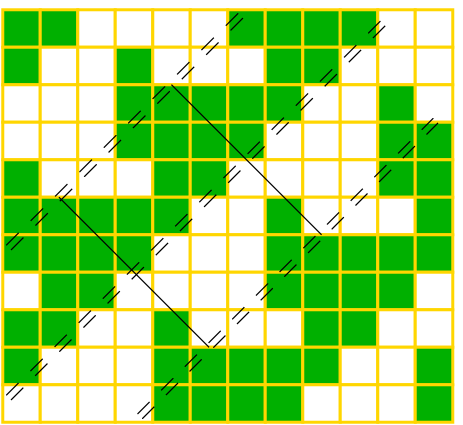}\hskip 10 pt\epsffile{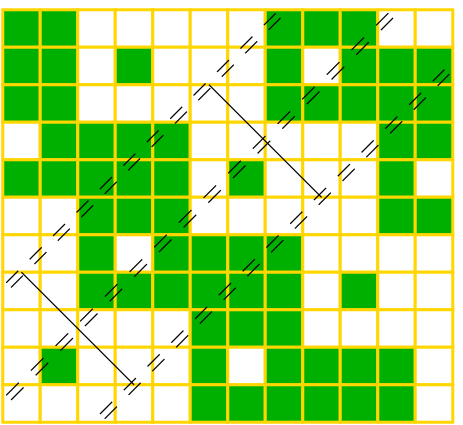}

\noindent (a)\hskip 1.8 in (b)

\noindent Figure 10. $G_1$ and $H_1$ lattice units. a. Order-24 example of species $1_e$ with $\ell =3$, $w=4$. \hskip 10 pt  b. Order-30 example of species $1_o$ with $\ell =5$, $w=3$.
\end{figure}

The production of examples such as those in Figure 10 is easy if it is possible.
To produce all the members of a family, all that is needed is to construct a diagram of the lattice unit and name enough cells that their (distinct) orbits under the group cover the plane.
This refines a technique introduced by Gr\"unbaum and Shephard [3].
In the case $\ell =3$, $w=4$, the family of Figure 10a, 12 cells/orbits $A$ to $L$ illustrated in Figure 11 are required.
\begin{figure}
\noindent
\epsffile{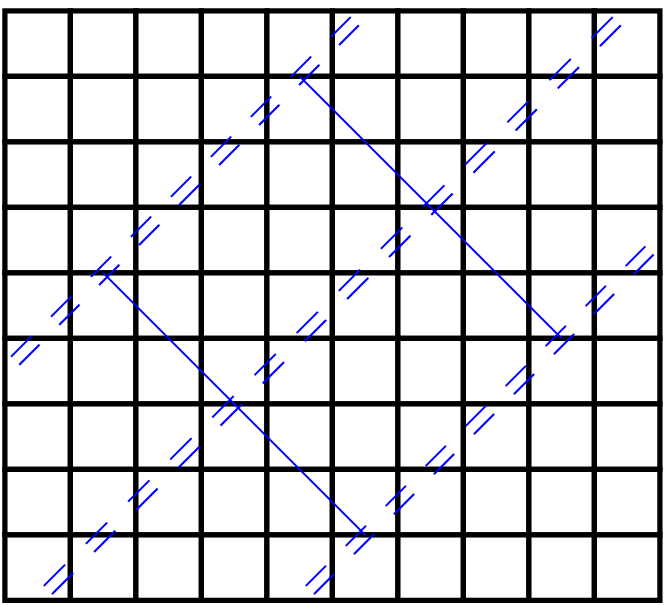}

\vskip -169 pt
\noindent
$\begin{array}{cccccccccc}
B & \overline{A} & \overline{D} & \overline{H} & \overline{K} & J & H & E & C & \overline{B}\\[6.7pt]
\overline{C} & \overline{F} & \overline{G} & \overline{J} & L & G & D & B & \overline{A} & \overline{D}\\[6.7pt]
\overline{E} & \overline{I} & \overline{L} & K & I & F & A & \overline{C} & \overline{F} & \overline{G}\\[6.7pt]
\overline{H} & \overline{K} & J & H & E & C & \overline{B} & \overline{E} & \overline{I} & \overline{L}\\[6.7pt]
\overline{J} & L & G & D & B & \overline{A} & \overline{D} & \overline{H} & \overline{K} & J\\[6.7pt]
K & I & F & A & \overline{C} & \overline{F} & \overline{G} & \overline{J} & L & G\\[6.7pt]
H & E & C & \overline{B} & \overline{E} & \overline{I} & \overline{L} & K & I & F\\[6.7pt]
D & B & \overline{A} & \overline{D} & \overline{H} & \overline{K} & J & H & E & C\\[6.7pt]
A & \overline{C} & \overline{F} & \overline{G} & \overline{J} & L & G & D & B & \overline{A}
\end{array}$

\medskip
\noindent Figure 11. The letter array for production of Figure 10a. ($X$ and $\overline{X}$ are of complementary colours.)
\end{figure}
The side-preserving (colour-reversing) glide-reflections make half of the plane be coloured with the complementary colours ${\overline A}$ to $\overline {L}$.
Colours have then to be chosen such that no further symmetry is introduced.
In Figure 10a, $H$ and $J$ are pale, other unbarred letters dark.
It is often possible to produce plain weave because its group contains many groups as subgroups.
When $\ell =3$, $w=2$, is attempted with six orbits, further symmetry is unavoidable.
This was confirmed by noting that, in the catalogue of Gr\"unbaum and Shephard [7], the only order-12 fabrics of species 1 are the two Roth gave as examples, 12-183-1 (Figure 4a) and 12-411-1.

\noindent {\it 2.}\hskip 8 pt When $G_1$ is of type $pg$ but illustrated by Figure 3d as having side-reversing (colour-preserving) glide-reflections so that $H_1$ has only translations ($p1$) and the Roth type is 2, little that was said of species 1 needs to be changed.
The possible $G_1$ lattice units are the same for the analogous species $2_m$, $2_e$, and $2_o$, including the failure to find an example of $2_e$ with $\ell=3, w=2$.
\begin{figure}
\noindent
\epsffile{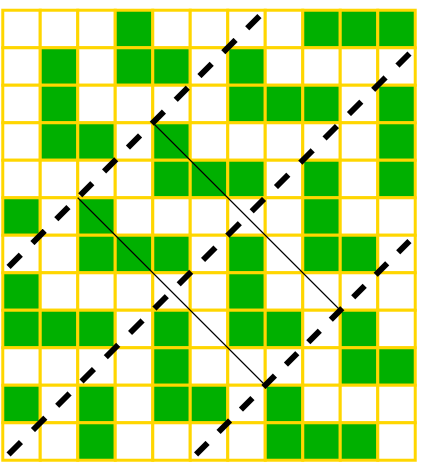}\hskip 10 pt\epsffile{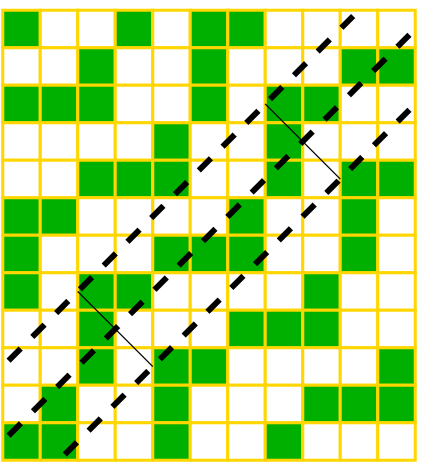}

\noindent (a)\hskip 1.65 in (b)

\noindent Figure 12. $G_1$ and $H_1$ lattice units. a. Order-20 example of species $2_m$ with $\ell=2$, $w=5$. \hskip 10 pt b. Order-20 example of species $2_e$ with $\ell=5$, $w=2$.
\end{figure}
Figure 12a shows a fabric of species $2_m$ ($\ell=2$, $w=5$) of smaller order than the example constructed by Roth ([1], Figure 3, $\ell=4$, $w=3$), and Figure 12b shows a fabric of species $2_e$ with $\ell=5$, $w=2$.
Period and order for all three species are $2\ell w$, for species $2_m$ and $2_e$ always divisible by 4 and for $2_o$ only by 2.
Genus is I.

\noindent {\it 3.}\hskip 8 pt When $G_1$ of type $pg$ is illustrated by Figure 3c as having a mixture of side-preser\-ving and side-reversing glide-reflections, with the latter absent from $H_1$, whose lattice unit is accordingly twice the width of $G_1$'s, we again have the glide-reflection axes with mirror position (type 3) or not (type 4).
In the former case, the lattice unit for $G_1$ must obey the constraints of species $1_m$.
One of Roth's two examples of species 3, which differ only insignificantly, is illustrated in Figure 8a and has in $\delta$ units the distance between neighbouring axes of the same kind $w=3$ and twice the length of the glides $\ell=2$.
So does 12-203-3, the only other example with so small a lattice unit.
The lattice unit of $H_1$, being twice as wide, has corners on only alternate strands, making these fabrics of pure genus II.
Period is $4\ell w$, but the order, on account of genus II, is only $2\ell w$, always divisible by 4.

\noindent {\it 4.}\hskip 8 pt When the axes of mixed side-preserving and side-reversing glide-reflections are not in mirror position, the $G_1$ lattice unit must obey the constraints of species $1_o$ or $1_e$, making the parameter values of $4_o$ and $4_e$ be those of $1_o$ and $1_m$ (rather than $1_e$ since $\ell=3$, $w=2$, is possible, being illustrated in Figure 8b). 
The $H_1$ lattice unit is twice as wide.
Examples of species $4_o$ with $\ell=5$, $w=3$, and $\ell=3$, $w=5$, are illustrated in Figure 13.
\begin{figure}
\noindent
\epsffile{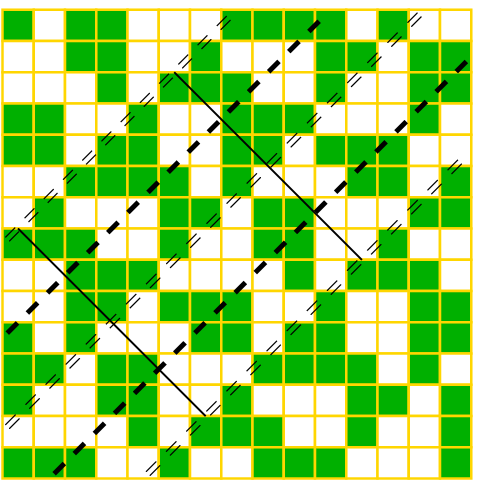}\hskip 10 pt\epsffile{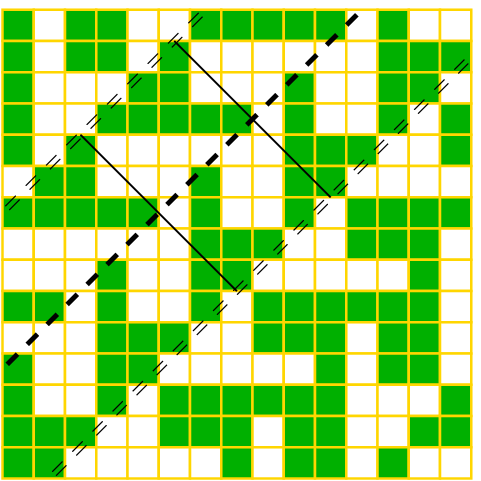}

\noindent (a)\hskip 1.86 in (b)
\smallskip

\noindent Figure 13. Order-60 examples of species $4_o$. a. $\ell=5$, $w=3$ ($H_1$). \hskip 10 pt b. $\ell=3$, $w=5$ ($G_1$ only).
\end{figure}
Too much symmetry results from $\ell =1$ or $w=1$.
The period and order of species $4_o$ are $4\ell w$, always divisible by 4, and of species $4_e$, with $\ell$ and $w$ from the parameters for $1_m$, also $4\ell w$, always divisible by 8.

\noindent {\it 5.}\hskip 8 pt We come now to groups of the Roth type of Figure 3e, species 5 with $G_1$ of type $pm$ and $H_1$ of type $p1$ since the side-reversing mirrors disappear from $H_1$. 
The mirrors are constrained as to position in cells but not otherwise, so that the isonemality constraint $(\ell ,w)=1$ is the only one.
Roth's examples (Figures 2a and 14a) show that $\{\ell,w\}=\{2, 3\}$ can be used both ways, and the same is true for the smallest viable odd pair, $\{3,5\}$ (Figures 14b and 14c).
\begin{figure}
\noindent
\epsffile{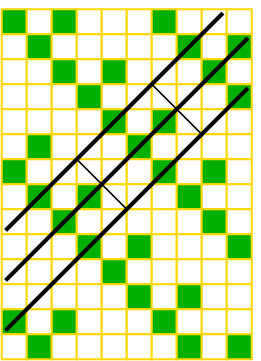}\hskip 10 pt\epsffile{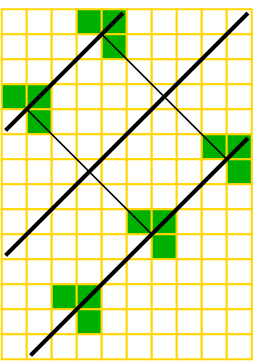}\hskip 10 pt\epsffile{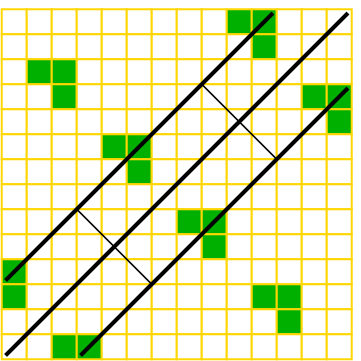}

\noindent (a)\hskip 0.99 in (b)\hskip 0.99 in (c)

\noindent Figure 14. $G_1$ and $H_1$ lattice units. a. 12-69-1 (species $5_e$, with $\ell=3$, $w=2$). \hskip 10 pt b, c. Order-30 examples of species $5_o$ with $\{\ell,w\}=\{3,5\}$, used both ways.
\end{figure}
The parameters for the obviously defined species $5_o$ and $5_e$ are those of $1_o$ and $1_m$ respectively.
While the period and order of species $5_o$ are $2\ell w$, always even, those of species $5_e$, also $2\ell w$, are divisible by 4.
Genus is I.

\noindent {\it 6.}\hskip 8 pt The second sort of type-$pm$ $G_1$ has axes of reflection that do not disappear from $H_1$ but become axes of (side-preserving) glide-reflection instead, making the $H_1$ lattice unit be twice as long as $G_1$'s and the symbol be $pm/pg$ (Figure 3f).
The constraints on these $G_1$ lattice units are those of type 5.
If the width of the lattice units is even, then the type is 6, illustrated in Figure 6.
The parameter values for type 6, in addition to those of $5_e$ (that is of $1_m$), are the pairing of 1 with even numbers, beginning with 2.
When the $G_1$ lattice-unit length, i.e., the length of the glide in $\delta$ units, $\ell=1$, prefabric 4-1-1* (the only prefabric of order 4 or less in species 1--10) has twice the distance between neighbouring mirrors $w=2$ and Roth's example 8-11-1 (Figure 6a) has $w=4$.
Period is $4\ell w$, but because the genus is II the order is only $2\ell w$, always divisible by 4.

\noindent {\it 7.}\hskip 8 pt On the other hand, if the width of the lattice units is odd $w\delta$, then the length $\ell\delta$ can be either odd or even giving rise to the two species $7_o$ and $7_e$ respectively.
The parameter values for species $7_o$, in addition to those of $5_o$ (that is of $1_o$), are the pairing of 1 with odd numbers, beginning with 3.
Roth's sole example illustrating type 7, 12-315-1, is of species $7_o$ with $\ell=1$, $w=3$, the only example with these smallest parameters.
Pairs with 1 can be used only that way.
The $1_o$ parameter pairs begin with $\{3, 5\}$. 
They can be used here both ways as is illustrated in Figure 15.
\begin{figure}
\noindent
\epsffile{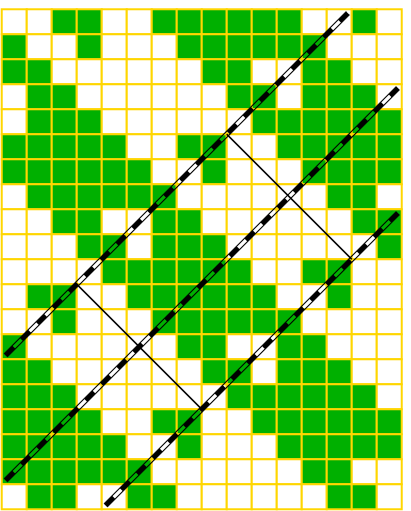}\hskip 10 pt\epsffile{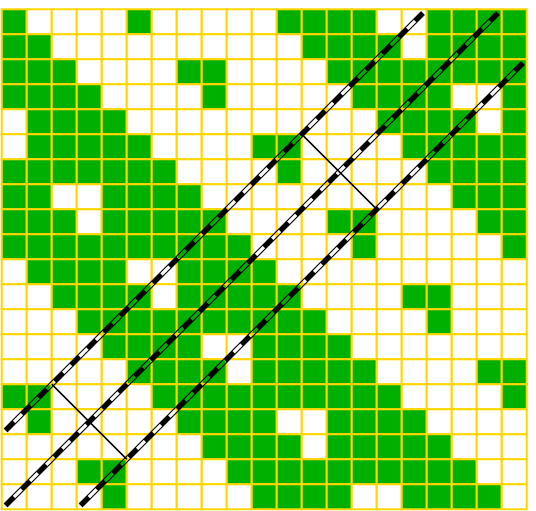}

\noindent (a)\hskip 1.57 in (b)

\noindent Figure 15. Order-60 examples of species $7_o$ with $\{\ell,w\}=\{3,5\}$, used both ways. $H_1$ lattice units.
\end{figure}
The parameter values for species $7_e$ are those of $5_e$ (that is of $1_m$), of which the smallest pair is $\ell = 2$, $w=3$, illustrated in Figure 16a.
\begin{figure}
\noindent\epsffile{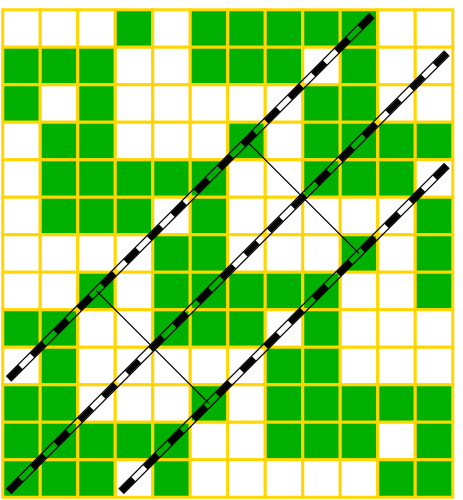}\hskip 10 pt\epsffile{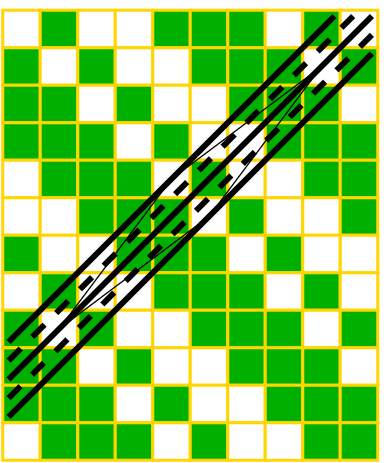}

\noindent Figure 16. $H_1$ lattice units.\hskip 10 pt a. Order-24 example of species $7_e$ with $\ell =2$, $w=3$.\hskip 10 pt b. The 1/1/2/3 twill, the smallest example of species $8_o$.
\end{figure}
Period and order for species $7_o$ are $4\ell w$, always divisible by 4, and for species $7_e$ are also $4\ell w$ but always divisible by 8 because $\ell$ and $w$ come from the $1_m$ parameter list.
The genus of both species $7_e$ and $7_o$ is I and II, as can be seen in Figure 15b. (The dark motif at the far right can go up one strand with colour reversal and down one strand without colour reversal by moving farther left.)
That completes the crystallographic type $pm$.

\noindent {\it 8.}\hskip 8 pt We turn now to $G_1$ of the type $cm$ illustrated in Figure 3g, in which mirrors and axes of side-reversing glide-reflections alternate.
The constraints on this group type for feasibility (parity) and isonemality (common factors) have already been discussed in the previous section, with Roth's examples illustrated in Figure 9.
The same parameters, greater than 1, can be used either way around.
$\{\ell, w\} = \{2, 4\}$ and $\{3, 5\}$, used one way in Figures 9a, b, are used the other way around in Figures 17a, b.
If the length $\ell\delta$ and width $w\delta$ of the rhombic lattice unit are both odd, then the parameters of what is species $8_o$, having only to be relatively prime, include those of species $1_o$, but a greater range of parameters is allowed than those of species $1_o$ since $w$ can be 1 so long as $\ell$ is at least 7. 
These extra parameters, like those of species 6 and $7_o$, can be used only that way around.
This adds, to what would otherwise be present, an infinite family of twills the smallest of which is the 1/1/2/3 twill of order 7 (Figure 16b).
If $\ell$ and $w$ are even, then the parameters of species $8_e$, their halves needing to be both relatively prime and of opposite parities, must be twice the parameters for type 6, that is twice those of species $1_m$ augmented by 1 paired with all even numbers.
Roth's examples (Figure 9) are of both kinds.
\begin{figure}
\noindent
\epsffile{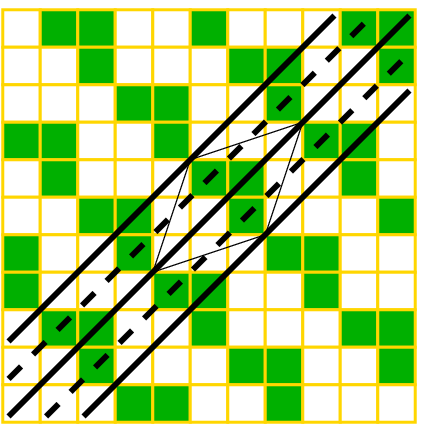}\hskip 10 pt\epsffile{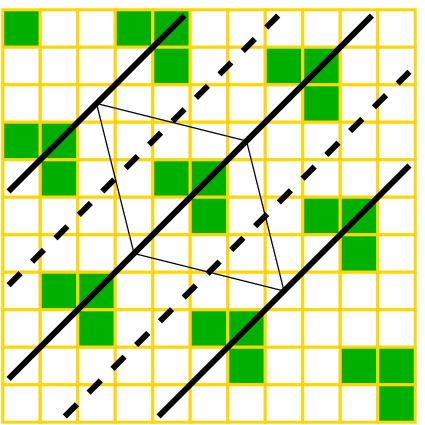}

\noindent (a)\hskip 1.65 in (b)

\noindent Figure 17. Examples of species $8_e$, $8_o$, with $G_1$ and $H_1$ lattice-unit dimensions of Figure 9 used the other way around.\hskip 10 pt a. 8-19-5. \hskip 10 pt b. 15-19-1.
\end{figure}
Period and order for species $8_o$ is $\ell w$, always odd, the only odd orders among these prefabrics.
For species $8_e$, period and order are also $\ell w$, always divisible by 8.
Genus is I.

\noindent {\it 9.}\hskip 8 pt Finally, $G_1$ of the type $cm$ illustrated in Figure 3h has $H_1$ of type $pg$ with lattice units either an even or odd multiple of $\delta$ wide.
If the width is an even multiple $w\delta$, then the Roth type is 9, the length of the rhomb $\ell\delta$ must be even too, and the axes of glide-reflection are in mirror position (Figure 7).
If the width and length are odd multiples of $\delta$, then the Roth type is 10, and the axes of glide-reflection are not in mirror position because of the distance between the mirrors (Figure 5).
The constraints on $G_1$ for species 9 are the same as those of species $8_e$, and so the parameters are the same.
Roth's examples (Figure 7) are of the smallest pair, $\{2,4\}$, used both ways.
While the parameters are the same as in $8_e$, the $H_1$ lattice unit is larger, being rectangular rather than rhombic (Figure 7b). 
Period is $2\ell w$, but, because genus is II, order is $\ell w$, always divisible by 8.

\noindent {\it 10.}\hskip 8 pt The constraints on $G_1$ for species 10 are the same as those of species $8_o$; again a greater range of parameters is allowed than those of species $1_o$ since the width can be $\delta$ so long as the length is at least $7\delta$. 
The infinite family of  twills added to what would otherwise be present has as a smallest example the 1/2/4/1/2/4 twill of order 14 illustrated in Figure 5 along with Roth's example.
Period and order are $2\ell w$, always even.
Genus is I and II.
\vs

\noindent {\bf 6. Doubling} 

\noindent A natural question that has, so far as I know, never been asked because there has been no way to look at it is which designs remain isonemal when doubled as plain weave is doubled to 4-3-1; each strand is replaced by a pair of strands with the same behaviour.
Plain weave remains isonemal when doubled because it has symmetries not being considered here.
Doubling is something that weavers actually do, but it was introduced into the weaving literature by Gr\"unbaum and Shephard [3] applied to square satins like 5-1-1.
(They used the term again with a different meaning in [4].)
Among the prefabrics of the species being considered here, the prevalence of a requirement of oddness makes it easy to see that replacing every strand with two strands behaving the same way---and thereby doubling lattice-unit dimensions---will alter the dimensions to destroy isonemality.
The only species that come close to working are $8_o$ and 10 since there are corresponding species $8_e$ and 9 with even dimensions. 
But those even dimensions are twice numbers of different parity, which the odd dimensions of $8_o$ and 10 are not.
\smallskip
\begin{Thm}{No prefabric with only parallel symmetry axes can be doubled and remain isonemal.}
\end{Thm}
\smallskip
\noindent {\it Proof.} If one thinks of each strand being divided into two lengthwise this becomes clear. 
For half a strand to be mapped to the other half of the strand by a symmetry as isonemality requires, the half-strand must be rotated (there are no half-turn symmetries in these prefabric species) or reflected (no axis parallel to strands is available) or translated.
While the other two possibilities are realized by plain weave, not even it allows translation of the lower half of a weft to the upper half of it by a symmetry of the design, for to have that symmetry every cell of the weft would have to be coloured identically to that above it since the weft's upper half would be mapped to the lower half of the next weft, cell images straddling the weft boundary.
The design would accordingly have to consist of vertical strips and so not be isonemal unless it were the trivial prefabric.
It is clearly impossible.
\vs

\noindent {\bf 7. Halving} 

\noindent By analogy with doubling we shall call {\it halving} the operation on a prefabric introduced by Gr\"unbaum and Shephard in [4], the removal of every other strand in each direction, making an intermediate construction used in chair seats (and the chapel floor at New College Oxford from 4-3-1) and called a pseudofabric by them, and the widening of all the strands uniformly to produce another prefabric.
Gr\"unbaum and Shephard left as an open question when the prefabric produced by halving is isonemal.
In terms of species, the answer for the prefabrics considered here is simple, but obtaining it will require a little machinery.
We require the cells of the plane to be numbered in fours: group the cells into square blocks of four that tessellate the plane as do the cells, and number the cells of each block as the quadrants of the Cartesian plane are numbered.
Plain weave can be described in these terms as odd dark and even pale or vice versa.
Beginning with a prefabric, one can then produce any half-fabric by preserving the crossings at cells assigned any one of the four numbers and discarding the rest.
If one keeps the 2-cells, for example, one has discarded the strands in both directions containing the 4-cells so that what were crossings in the odd cells are no longer crossings because what was crossed by the remaining strands is no longer there.
And the areas of the plane occupied by the odd cells and the 4-cells have been shrunk to nothing by widening the remaining strands to expand the 2-cells to fill the plane.
Since there is a sense in which the prefabric obtained by keeping the 2-cells and that obtained by keeping the 4-cells are combined (by the crossings in the odd cells) to make the original prefabric, I shall refer to the results of this process as {\it factors} and the 2-factors and the 4-factors collectively as the {\it even} factors and the others as the {\it odd} factors.
The action may all be in the combination: if halving is applied to plain weave, all four factors are trivial prefabrics.
This is, in fact, a hint of one result, namely that if the $G_1$ lattice unit of the prefabrics under consideration here is small enough, a twill (of which the trivial prefabric is the trivial example) results from halving.
Small enough is a length or width of $2\delta$ because halved it is $\delta$, and that forces twills, all of which are isonemal.

All prefabrics of species 1 to 10 are of genus I or II.
Accordingly halved they become of genus I with translations from weft to adjacent weft.
Moreover, the double offset from a strand to a strand {\it two} strands away is always even so that the numbering of the cells is preserved under those translations.
So in a halved prefabric, there are always strand-to-adjacent-strand translations that are symmetries. 
The question is whether there are warp-to-weft symmetries.
The species here all have glide-reflections or reflections that take warp to weft, and so the question reduces to whether they preserve the numbering.
Glide-reflections with axes not in mirror position interchange odd and even cells, and so they cannot be called upon to make a halved prefabric isonemal.
Prefabrics with only such warp-to-weft symmetries (species $1_o$, $1_e$, $2_o$, $2_e$, and 4) will in consequence not halve to isonemal factors unless their lattice units are small enough.
Ignoring those species then, we need to see what reflections and mirror-position glide-reflections do to the numbering of the cells.

\begin{table}
\begin{tabbing}
\noindent\=\hskip 0.75 in  \= Parities\hskip 10 pt \=\hskip 0.4 in \=\hskip 1 in\=\hskip 0.4 in\=\\
\>\> of axis\>\> Same\>\>Opposite\\
\>\> and cells\\
\>Operation\\
\>Odd glide \>\> interchanges numbers\>\>preserves numbers\\
\>Even glide \>\> preserves numbering\>\>interchanges numbers\\
\>Reflection \>\> preserves numbering\>\>interchanges numbers\\
\>(even length)\\
\>Reflection \>\> preserves numbering\>\> preserves numbering when\\
\>(odd length)\>\>\>\> combined with translation
\end{tabbing}

\smallskip
\noindent Table 1. Effects of reflections and glide-reflections on cell numbers.
\end{table}

\begin{Lem}{The action of a reflection or glide-reflection with axis in mirror position on the numbering of the cells is displayed in Table 1, where the {\it parity} of an axis is the parity of the cells through which it runs, {\it length} refers to length of the $G_1$ lattice unit, and {\it interchange of numbers} means specifically the interchange of $1$ and $3$ or $2$ and $4$.}
\end{Lem}
\smallskip

\noindent {\it Proof.} The effects of glide-reflections in Table 1 are obvious. 
And the effects of reflections are those of zero (even) glides except for the anomalous effect on cells of parity opposite to that of the cells through which the mirror runs when the $G_1$ lattice unit is of odd length. 
Then, while the reflection interchanges numbers, a translation along the mirror by the length of the lattice unit (a symmetry of the prefabric perhaps involving $\tau$) also interchanges numbers (e.g., Figure 17b). 
Since the product of the two interchanges is the identity, the numbering is preserved by a subgroup of the original symmetry group.

\smallskip
\begin{Thm}{Isonemal prefabrics of species $1$ to $10$, except those satisfying the two conditions (1) having only glide-reflections with axes not in mirror position (species $1_o$, $1_e$, $2_o$, $2_e$, and $4$) and (2) having both lattice-unit dimensions greater than $2\delta$, can be halved with only isonemal results.}
\end{Thm}

\smallskip
\noindent {\it Proof.} 
In species 1 to 10 (except for $1_o$, $1_e$, $2_o$, $2_e$, and $4$) there are four different situations with respect to the spacing of axes:

\noindent 1. There are only glide-reflection axes, and they are half an odd multiple of $\delta$ apart (species $1_m$, $2_m$, and 3).

\noindent 2. There are mirrors, and they are half an odd multiple of $\delta$ apart (species $5_o$, part of $5_e$ with $\ell$ even, 7, $8_o$, 10).

\noindent 3. There are mirrors a whole multiple of $\delta$ apart and with $G_1$ lattice units an odd multiple of $\delta$ long (the part of species $5_e$ with $\ell$ odd and species 6).

\noindent 4. There are mirrors an even number of $\delta$ apart with glide-reflection axes between them and with the glides an odd multiple of $\delta$ (species $8_e$ and 9).

In all cases, the result of halving is isonemal. What needs to be proved is the presence of a warp-to-weft transformation; case by case:

\noindent 1. Whether the glide-reflection axes are for odd glides or even, adjacent axes run through even cells and odd cells, and so factors of both parities have their numbering preserved by one or other of any pair of adjacent axes' glide-reflections by Lemma 5.

\noindent 2. Adjacent mirrors run through even cells and odd cells, and so factors of both parities have their numbering preserved by reflection in one or other of any pair of adjacent mirrors by Lemma 5. (The glide-reflections of species $8_o$ and 10, whose axes are not in mirror position, take no part in this and vanish from the factors.)

\noindent 3. While, because of their spacing, all mirrors run through cells of factors of the same parity, preserving the numbering of cells in those factors, the effect of the mirrors on the factors of the other parity is that of glide-reflections by Lemma 5. (Since the cells that the mirrors run through are not in those factors, the axes of glide-reflection are not in mirror position in those factors. They are side-preserving or side-reversing depending on whether $\tau$ is involved in the translation from one $G_1$ lattice unit to the next in the direction of the axes.)

\noindent 4. Both mirrors and axes of glide-reflection run through cells of factors of the same parity. The mirrors preserve the numbering of the factors they run through, and the glide-reflections, their glides being odd, preserve the numbering of the other factors by Lemma 5.

The theorem is proved.
\vs

\noindent {\bf 8. Fabrics of a given order} 

\noindent In the absence of Roth's classification of symmetry-group types and the above determination of the possible symmetry groups, isonemal fabrics of various orders had to be found by a brute-force method relying on trial and error. 
This was explained at the beginning of Gr\"unbaum and Shephard's catalogue [6].
While period is perhaps a more natural way to classify fabrics than order, it is now easy to find all fabrics of a given period or order with no trial and error, only the elimination of fabrics of smaller (divisor) period or order and prefabrics that fall apart, both of which can be determined.
This is at present possible only for the species discussed here, but this discussion will be extended.

An odd order is a prime or a prime power (and so has no relatively prime factors) or else it does have relatively prime factors.
In the first and second cases, numbers like 23 and 25, there are only the twills of those orders in species $8_o$ available.
In the third case, for a number like 21, one looks {\it as well} for its factorizations among the parameters for the only species allowing non-twills of odd order, $8_o$, namely the parameters for $1_o$.
Since they are only 3 and 7, used both ways, each possible lattice unit and therefore its unique family are determined.
The method for finding the members of a family explained in section 5 can then be used to determine all of the fabrics in the family.
Since these factors have no factors, no extraneous fabrics or prefabrics will be produced.
When a factor has factors, whole families of smaller-order designs need to be rejected.

For an even order $N$, one first looks for the species that allow that order by looking for the possibilities of $N=\ell w$, $2\ell w$, or $4\ell w$ in all appropriate rows of Table 2 (combining information stated above in the discussion of each species).
Species 6, $7_o$, and 10 add combinations of 1 with even, odd, and odd numbers, respectively, provided they are 2 or more, 3 or more, and 7 or more, respectively to the parameters of $1_m$, $1_o$, and $1_o$, respectively.
The plethora of possibilities explains why the catalogue stopped at orders 15 and 17; even without trial and error, the number of possibilities of order 16 is too large.

\begin{table}
\begin{tabbing}
Condition\hskip 10 pt \= Order\hskip 10 pt \= Parameters\hskip 10 pt \= Species\\
$2 | N$, $4\!\not |\, N$ \> $2\ell w$ \> $1_o$         \> $1_o, 2_o, 5_o, 10$\\
                         \> $2\ell w$ \> 10            \> 10\\
$4 | N$, $8\!\not |\, N$ \> $4\ell w$ \> $1_o$         \> $4_o, 7_o$\\
                         \> $4\ell w$ \> $7_o$         \> $7_o$\\
$4 | N$                  \> $2\ell w$ \> $1_m$         \> $1_m, 1_e, 2_m, 2_e, 3, 5_e, 6$\\
                         \> $2\ell w$ \> $6$            \> 6\\
$8 | N$                  \> $\ell w$  \> $8_e (1_m, 6)$ \> $8_e$, 9\\
                         \> $4\ell w$ \> $1_m$          \> $4_e, 7_e$
\end{tabbing}

\noindent Table 2. Species in which prefabrics of even orders $N$ may be found. 
\end{table}

If one runs through the orders from 5 to 27, one finds the following results:

\noindent 5 is too small to use for the only possibility, with 1 for a twill of species $8_o$.

\noindent 7, 11, 13, 17, 19, and 23 are odd primes, and 9, 25, and 27 are powers of primes and large enough to be used for species-$8_o$ twills.

\noindent $6=2(1\times 3)$, and $10=2(1\times 5)$ are too small for species 10.

\noindent $14=2(1\times 7)$, $18=2(1\times 9)$, $22=2(1\times 11)$, and $26=2(1\times 13)$ are big enough that they can be used for species 10.

\noindent $15= 3\times 5$ and $21=3\times 7$ can be used in factored form for species $8_o$ both ways around, and each can be used with 1 for a species-$8_o$ twill.

\noindent $8=2(1\times 4)$, $12=2(1\times 6)$, $16=2(1\times 8)$, $20=2(1\times 10)$, and $22=2(1\times 11)$ can be used for species 6 one way.
The designs of order 8 will all reappear among the designs of order 16 if not eliminated since $G_1$ of these order-16 designs is a subgroup of the $G_1$ of these order-8 designs.
The nesting of the symmetry groups, which happens again and again, is an unavoidable geometrical fact; it is not reflected in the way I have chosen to name the species.

\noindent $8=2\times 4$, $16=2\times 8$, and $24= 2\times 12=4\times 6$ can be used for species $8_e$ and 9 both ways around.

\noindent $12=2(2\times 3)$, $20=2(2\times 5)$, and $24=2(4\times 3)$ can be used for species $1_m, 1_e, 2_m, 2_e,$ $3, 5_e$ and 6.

\noindent $12=4(1\times 3)$ and $20=4(1\times 5)$ can be used for species $7_o$ one way.

\noindent $24= 4(2\times 3)$ can be used for species $4_e$ and $7_e$ one way.

\noindent One notices also that, because of what the parameters need to be, species $1_o, 2_o, 5_o$ and the 10s with the $1_o$ parameters begin with order $2\ell w=30$, and that species $4_o$ and the $7_o$s with the $1_o$ parameters begin with order $4\ell w=60$ and so are not represented among the prefabrics of orders 5--27.

Having determined which species can contain a family of the desired order, one can determine all of the prefabrics in each such family by the method explained in section 5.
There will typically arise in such determination three classes of design that are not wanted.
There may be designs with order equal to a factor of the desired order.
There may be designs with the symmetry group chosen being a subgroup of their respective groups, of which the previous sentence mentions a subset.
And there can be prefabrics that fall apart in the cases, pointed out by Roth [1], of species 3, 6, and 9, in all of which the order is divisible by 4.

One can take care to avoid constructing sets of orbits that have additional symmetry, or one can construct all that there appear to be and reject those with additional symmetries.

One can take care to avoid constructing prefabrics that fall apart, or one can construct all that there appear to be and reject those found to fall apart.
Or one can construct them beforehand and avoid or reject them as they might arise.
Species 3 has even $\ell$ beginning with 2 and odd $w$ beginning with 3. 
Possible orders of prefabrics that fall apart are $2\ell w$, beginning with 12.
The initial example is 12-69-2*, there being no other of order less than 20.
Species 6 has odd $\ell$ beginning with 1 and even $w$ beginning with 2, but the smallest pair for which there is a prefabric that falls apart has $\ell=1$, $w=6$.
Order is $2\ell w$.
The initial examples are 12-69-1* and 16-1093-1* with $\ell=1$ and $w=6$ and 8 respectively.
There are no others of order less than 20.
Species 9 has even $\ell$ and $w$ and order $\ell w$ divisible by 8, the smallest feasible order. 
The smallest actual prefabrics of this species that fall apart are 16-277-4* and 16-1093-3* with $\ell=2$ and $w=8$.
There are no others of order less than 24.

The prefabrics that fall apart are easily generated by the usual method, noting that alternate rows have every other cell dark and every other cell pale with these predetermined cells forming a checkerboard with the undetermined cells.
The mirrors and axes of glide-reflection can be placed onto this partially determined array only consistently with its colouring.
Along diagonals of cells with predetermined colouring can be placed mirrors, axes of side-preserving glide-reflection with odd glides, and axes of side-reversing glide-reflection with even glides.
Between those diagonals, that is, through cells not predetermined, can be placed axes of side-preserving glide-reflection with even glides and axes of side-reversing glide-reflection with odd glides.
The parameters of these three species make these choices consistent.
The same generation of unwanted designs with more symmetry than wanted occurs among prefabrics that fall apart; for example, 12-21-2* with symmetry of Roth type 23 arises when one aims at order 12 and species 3 because of unwanted mirrors perpendicular to the imposed axes of glide-reflection.

The falling apart of such prefabrics into two components is a species of halving.
The warps corresponding to predominantly dark columns of the design---together with the wefts corresponding to predominantly pale rows---will lift off the other half of the warps and wefts.

The results of sections 6--8 need to be extended to the remaining prefabrics.
\vs

\noindent {\bf References}

\noindent [1] Richard L.~Roth, The symmetry groups of periodic isonemal fabrics, Geom.~Dedicata 48 (1993), 191-210.

\noindent [2] ------ Perfect colorings of isonemal fabrics using two colors, Geom.~Dedicata 56 (1995), 307-326.

\noindent [3] Branko Gr\"unbaum, Geoffrey C.~Shephard, Satins and twills: An introduction to the geometry of fabrics, Math.~Magazine 53 (1980), 139-161.

\noindent [4] ------ Isonemal fabrics, Amer.~Math.~Monthly 95 (1988), 5-30.

\noindent [5] D.~Schattschneider, The plane symmetry groups: Their recognition and notation, Amer.~Math.~Monthly 85 (1978), 439-450.

\noindent [6] Branko Gr\"unbaum, Geoffrey C.~Shephard, A catalogue of isonemal fabrics, in: Jacob E.~Goodman {\it et al.} (Eds.), Discrete Geometry and Convexity, Annals of the New York Academy of Sciences 440 (1985), 279-298.

\noindent [7] ------ An extension to the catalogue of isonemal fabrics, Discrete Math.~60 (1986), 155-192.

\noindent [8] C.R.J.~Clapham, When a fabric hangs together, Bull.~London Math.~Soc.~12 (1980), 161--164.

\noindent [9] T.C.~Enns, An efficient algorithm determining when a fabric hangs together, Geom. Dedicata 15 (1984), 259--260.

\noindent [10] W.D.~Hoskins, R.S.D.~Thomas, Conditions for isonemal arrans on a cartesian grid, Linear Algebra Appl. 57 (1984), 87--103.

\noindent [11] Cathy Delaney, When a fabric hangs together, Ars Combinatoria 21A (1986), 71--79.

\noindent [12] J.A.~Hoskins, R.S.D.~Thomas, The patterns of the isonemal two-colour two-way two-fold fabrics, Bull.~Austral.~Math.~Soc.~44 (1991), 33-43.

\noindent [13] Branko Gr\"unbaum, Geoffrey C.~Shephard, Tilings and patterns, W.~H.~Freeman, New York, 1987.

\end{document}